\documentclass[11pt]{amsart}
\usepackage{amsmath, amssymb, amsthm, amscd, graphicx}
\usepackage{epstopdf}

\theoremstyle{plain}
\newtheorem{prop}{Proposition}[section]

\newtheorem{lemm}[prop]{Lemma}

\theoremstyle{definition}
\newtheorem{defi}[prop]{Definition}

\newtheorem{nota}[prop]{Notation}
\newtheorem{exam}[prop]{Example}
\newtheorem{rema}[prop]{Remark}

\numberwithin{equation}{section}

\def\Reff#1; #2; #3; #4; #5; #6; #7\par{%
\bibitem{#1} #2, {\it #3}, #4 {\bf #5} (#6) #7}

\def\Ref#1; #2; #3; #4\par{%
\bibitem{#1} #2, {\it #3}, #4}

\catcode`\¨=\active\def¨{~}

\renewcommand{\aa}[1]{a_{#1}}
\newcommand{\ALD}{A\!L\!D}

\newcommand{\Bb}{B_\bullet}
\newcommand{\Bi}{B_\infty}
\newcommand{\bx}{\beta}
\newcommand{\by}{\gamma}

\newcommand{\concat}{{}^{\scriptscriptstyle\frown}}

\newcommand{\dd}{...}
\newcommand{\ddd}{, ...\,,}

\newcommand{\eALD}{=_{\scriptscriptstyle\!A\!L\!D}}
\newcommand{\eLD}{=_{\scriptscriptstyle\!L\!D}}
\newcommand{\eLDs}{\eLD}
\newcommand{\ex}{x}
\newcommand{\ey}{y}
\newcommand{\ez}{z}

\let\ge=\geqslant
\renewcommand{\gg}{\gamma}

\newcommand{\hR}{\mathrm{ht}_{\!\scriptscriptstyle R}}

\newcommand{\ie}{{\it i.e.}}
\newcommand{\II}[1]{I(#1)}
\renewcommand{\Im}{\mathrm{Im}}
\newcommand{\inv}{^{-1}}

\newcommand{\JJ}[1]{J(#1)}

\newcommand{\LD}{L\!D}
\let\le=\leqslant
\renewcommand{\lg}[1]{\ell(#1)}

\newcommand{\oop}{\cdot}
\newcommand{\op}{*}
\newcommand{\Op}{\mathbin{\scriptstyle\square}}
\newcommand{\OP}{\circ}
\newcommand{\ops}{\mathbin{\vec\op}}

\newcommand{\pLD}{\pref_{\scriptscriptstyle L\!D}}
\newcommand{\pref}{\mathrel{\scriptstyle\sqsubset}}
\newcommand{\prefs}{\mathrel{\vec{\pref}}}

\newcommand{\resp}{{\it resp.{¨}}}

\let\s=\sigma
\newcommand{\Seq}[1]{\widehat{#1}}
\newcommand{\sh}{\partial}
\newcommand{\ssh}{\partial}
\newcommand{\sm}{<}
\newcommand{\smALD}{<_{{}_{\!A\!L\!D}}}
\renewcommand{\ss}[1]{\sigma_{#1}}
\renewcommand{\SS}{S}

\newcommand{\To}[1]{T^{\op}_{#1}}
\newcommand{\TO}[1]{T^{\OP}_{#1}}
\newcommand{\TOo}[1]{T^{\op,\OP}_{#1}}
\newcommand{\ts}{s}
\renewcommand{\tt}{t}
\newcommand{\tti}{t_1}
\newcommand{\ttii}{t_2}
\newcommand{\ttiii}{t_3}
\newcommand{\TTo}{T^{\op}}
\newcommand{\TTO}{T^{\OP}}
\newcommand{\TTOo}{T^{\op,\OP}}
\newcommand{\tu}{u}
\newcommand{\tv}{v}

\newcommand{\vts}{\,\,\vec{\vrule width0pt height5pt}\!\ts\,}
\newcommand{\vtt}{\,\,\vec{\vrule width0pt height5pt }\!\tt\,}

\newcommand{\xx}{x}
\newcommand{\xs}{s}
\newcommand{\xt}{t}
\newcommand{\xu}{u}

\newcommand{\yy}{y}

\newcommand{\zz}{z}
\begin{document}$$$$

\author{Patrick DEHORNOY}
\address{Laboratoire de Math\'ematiques Nicolas
Oresme UMR 6139\\ Universit\'e de Caen,
14032~Caen, France}
\email{dehornoy@math.unicaen.fr}
\urladdr{//www.math.unicaen.fr/\textasciitilde dehornoy}

\title{Free augmented LD-systems}



\begin{abstract}
Define an augmented LD-system, or ALD-system, to be a set
equipped with two  binary operations, one satisfying the left
self-distributivity law $x \op (y \op z) = (x \op y) \op (x \op z)$
and the other satisfying the mixed laws $(x \OP y)
\op z = x \op (y \op z)$ and $x \op (y \OP z) = (x \op y) \OP (x
\op z)$. We solve the word problem of the ALD laws, and prove
that every element in the parenthesized braid group~$\Bb$
of~\cite{Bri1, Bri2, Dhb, Dhe} generates a free ALD-system of
rank¨$1$, thus getting a concrete realization of the latter
structure.
\end{abstract}

\maketitle

Define an {\it LD-system} to be an algebraic system
made of a set~$\SS$ equipped with a binary
operation~$\op$ that satisfies the left self-distributivity
law
\begin{equation}
\tag{$\LD$} x \op (y \op z) = (x \op y) \op (x \op z).
\end{equation}
Classical examples include groups equipped with their
conjugacy operation $\xx \op \yy = \xx \yy
\xx\inv$, and lattices with their inf or sup
operation. Less classical examples have appeared in Set
Theory with the iterations of elementary embeddings
\cite{Lvb}, and in Low Dimensional Topology where $(\LD)$
provides an algebraic translation of Reidemeister move~III
\cite{Joy, Mat, FeR}. A rich theory has been developed for
LD-systems~\cite{Dgd}. In particular, it is known that
there exists on Artin's braid group~$\Bi$ an
LD-operation~$\op$ such that the $\op$-closure of any
braid is a free LD-system of rank~$1$---which provides
a concrete realization of the latter structure.

Many examples of LD-systems turn out to be equipped
with a second operation connected in various ways with
the self-distributive operation. In the typical case of
group conjugacy, using~$\circ$ for the group product, the
following mixed identities are satisfied
\begin{gather}
\tag{$\ALD_1$}
\xx \op (\yy \op \zz) = (\xx \OP \yy)
\op \zz,\\
\tag{$\ALD_2$} \xx \op (\yy \OP \zz) = (\xx \op \yy)
\OP (\xx \op \zz).
\end{gather} 
When we add the identity $\xx \OP \yy = (\xx \op \yy) \OP
\xx$, the associativity of~$\circ$ and the existence of
a unit, one obtains the structure of an LD-monoid,
which is investigated in Chapter~XI of~\cite{Dgd} (and
in¨\cite{Drc, Drd} under the name of LD-algebra).

It is easy to verify that all LD-systems cannot be enriched into
LD-monoids. In particular, this is the case for the above
mentioned LD-structure on~$\Bi$, for which there can exist no
second operation verifying¨$(\ALD_1)$. In~\cite{Dhe}, building
on earlier approaches of~\cite{Bri1, Bri2, Dhb}, a new
group~$\Bb$ extending both Artin's braid group~$\Bi$ and
R.Thompson's group~$F$ is investigated. This group is called the
parenthesized braid group, as its elements can be
naturally interpreted using braid diagrams in which the
strands come grouped into blocks that can be encoded
in parenthesized words. It is shown that the
LD-structure of~$\Bi$ extends to~$\Bb$ and
that the latter can be completed with a second
operation that satisfies the above identities~$(\ALD_1)$
and~$(\ALD_2)$---but none of the further laws defining an
LD-monoid. Such a structure is called an {\it augmented}
LD-system, or {\it ALD-system}.

The aim of this note is to prove two new results about
ALD-systems: firstly, we solve the associated word problem,
and, secondly, we prove for the ALD-structure of the parenthesied
braid group~$\Bb$ a result similar to that established
in~\cite{Dgd} for the LD-structure of ordinary braids, namely
that every element of~$\Bb$ generates a free ALD-subsystem
of~$\Bb$. Being quite similar to those holding for~$\LD$
and~$\Bi$, these results are not surprising. However, their proofs
require a few new specific arguments that are the subject of this
paper.

\section{Free augmented LD-systems}

The aim of this section is to solve the word problem for the ALD
laws, \ie, to describe an algorithm that enables one to decide
whether two terms are or not equivalent up to¨ALD.

\subsection{ALD-systems}

The algebraic systems considered here are as follows:

\begin{defi}
An {\it ALD-system} is defined to be a set~$\SS$
equipped with two binary operations, $\op$
and~$\circ$ that satisfy the identities $(\LD)$, $(\ALD_1)$,
and¨$(\ALD_2)$.
\end{defi}

\begin{exam} \label{X:ALD}
We already observed that any group~$G$ equipped with
the conjugation operation~$\op$ and the product is an
ALD-system---and even an LD-monoid. Another easy
example is obtained by starting with an arbitrary
binary system $(\SS, \OP)$ and considering an
$\OP$-endomorphism~$f$. Then defining $\xx \op
\yy = f(\yy)$ turns $(\SS, \op, \OP)$ into an
ALD-system. 
\end{exam}

If $L_\xx$ denotes the left $\op$-translation
$\yy \mapsto \xx \op \yy$, then  $(\LD)$ and¨$(\ALD_2)$
express that, for each¨$\xx$ in the considered domain, $L_\xx$
is an endomorphism with respect to¨$\op$ and¨$\OP$,
respectively, while $(\ALD_1)$ expresses that
$\circ$ corresponds to a composition of translations:
$L_{\xx \OP \yy} = L_\xx \OP L_\yy$. Thus, an ALD-system is an
LD-system where the family of left translations is closed under
composition---and in which $(\ALD_2)$ is satisfied. It may be
noted that, in any case, the conjunction of~$(\LD)$
and~$(\ALD_1)$ implies some weak form
of~$(\ALD_2)$, as we can write
\begin{multline*}
(\xx \op (\yy \OP \zz)) \op (\xx \op \tu) 
=_{\LD} \xx \op ((\yy \OP \zz) \op \tu)
=_{\ALD_1} \xx \op (\yy \op (\zz \op \tu))\\
=_{\LD} (\xx \op \yy) \op ((\xx \op \zz) \op (\xx \op
\tu)))
=_{\ALD_1} ((\xx \op \yy) \OP (\xx \op \zz)) \op
(\xx \op \tu),
\end{multline*}
which follows from~$(\ALD_2)$ and actually implies it
if we may cancel $\xx \op \tu$ on the right.

\subsection{Terms and free ALD-systems}

We consider in the sequel free ALD-systems. As usual,
the latter can be introduced as quotients of absolutely free
algebras, \ie, of algebras consisting of terms subject to no
relation. Our notation will be as follows. 

\begin{defi}
For $n \ge 1$, we denote by~$\To n$ (\resp $\TO n$, \resp
$\TOo n$) the set of all binary terms constructed using the
operator~$\op$ (\resp $\OP$, \resp $\op$ and~$\OP$) from
$n$~fixed variables $\xx_1\ddd\xx_n$. We write
$\TTo$ for the union of all~$\To n$, and
similarly with~$\TTO$ and~$\TTOo$---and¨$\xx$
for¨$\xx_1$.
\end{defi}

The {\it size} of a term¨$\tt$ is defined to be the number of
occurrences of variables in¨$\tt$, \ie, it is defined to be¨$1$
when $\tt$ is a variable, and to be the sum of the sizes of the left
and the right subterms of¨$\tt$ otherwise. By construction,
$\TTOo n$ is an absolutely free algebra of rank~$n$.
The following is clear:

\begin{lemm} \label{L:Free}
Let $\eALD$ be the congruence on~$\TOo n$
generated by all instances of the laws¨$(\LD)$,
$(\ALD_1)$, and¨$(\ALD_2)$\footnote{\ie, all pairs of terms of
the form $(\tti \op (\ttii \op \ttiii), (\tti \op \ttii) \op (\tti \op
\ttiii))$,  $(\tti \op (\ttii \op \ttiii), (\tti \OP \ttii) \op
\ttiii)$, and $(\tti \op (\ttii \OP \ttiii), (\tti \op
\ttii) \OP (\tti \op \ttiii))$}. Then, for
each~$n$, the system $\TOo n/\!\!\eALD$ is a
free ALD-system of rank~$n$.
\end{lemm}

We say that two terms¨$\tt, \tt'$ are {\it ALD-equivalent} if $\tt
\eALD \tt'$ holds. Of course, there is a similar result for the
free LD-system of rank¨$n$ obtained as $\To n/\!\!\eLD$, where
$\eLD$ is the congruence generated by the instances of the sole
law~$(\LD)$.

It is helpful for intuition to associate with every term a finite
binary rooted, labeled tree: the tree associated with
a variable~$\xx$ consists of a single node labeled~$\xx$; for
$\Op = \op$ or~$\OP$, the tree associated with $\tti \Op \ttii$
consists of a root labeled~$\Op$ admitting as its left
subtree the tree associated with~$\tti$, and as its right
subtree the tree associated with~$\ttii$. 

As a preliminary remark, let us observe that the variety
of ALD-systems is properly intermediate between
LD-systems and LD-monoids.

\begin{prop}
$(i)$ A free LD-system cannot be enriched
into an ALD-system.

$(ii)$ A free ALD-system does not obey the
law~$\xx \OP \yy = (\xx \op \yy) \OP \xx$, and therefore is not
an LD-monoid.
\end{prop}

\begin{proof}
$(i)$ For~$\tt$ a term, let~$\hR(\tt)$ be the length of
the rightmost branch in the associated tree, \ie, define
$\hR(\tt)$ by $\hR(\xx) = 0$ and $\hR(\tti \Op \ttii) =
\hR(\ttii) + 1$ for $\Op = \op$ or¨$\OP$. Then the law¨$(\LD)$
preserves~$\hR$, and, therefore, $\hR$ induces a well defined
parameter on each free LD-system. On the other hand,
$(\ALD_1)$ changes¨$\hR$, so there may exist no
operation~$\OP$ satisfying~$(\ALD_1)$ on a free LD-system.

$(ii)$ The terms $\xx_1 \OP \xx_2$ and $(\xx_1 \op \xx_2) \OP
\xx_1$ are not ALD-equivalent, as none of the identities $(\LD)$,
$(\ALD_1)$, $(\ALD_2)$ may apply to a term with only two
occurrences of variables.
\end{proof}

\subsection{Two ALD-invariants}

In order to subsequently solve the word problem
of¨ALD, we shall associate with every term
in~$\TTOo$ two ALD-invariants, \ie, two objects
that depend only on the ALD-class of the term. The
first invariant is a term in~$\TO1$; the second one is
a finite sequence of LD-classes of terms in~$\TTo$.
To introduce the latter, we first fix some notation for
sequences.

\begin{nota}
Assume that $(\SS, \op)$ is a binary system. The set of
all finite, nonempty sequences of elements of~$\SS$ is
denoted by~$\Seq\SS$. An element of~$\Seq\SS$ is
typically denoted~$\vec\xs$; its length is then
denoted~$\lg{\vec\xs}$, and its successive elements
$\xs_1\ddd\xs_{\lg{\vec\xs}}$. The concatenation of
two sequences~$\vec\xs$, $\vec\xt$, \ie, the sequence of
length $\lg{\vec\xs} + \lg{\vec\xt}$ obtained by
writing~$\vec\xt$ after~$\vec\xs$, is denoted $\vec\xs
\concat \vec\xt$. Next, we denote
by~$\ops$ the binary operation on~$\Seq\SS$ defined
by
\begin{equation}
\vec\xs \ops \vec\xt = (\xs_1 \op \dd \op \xs_{\lg{\vec\xs}}
\op \xt_1 \ddd \xs_1 \op \dd \op \xs_{\lg{\vec\xs}}
\op \xt_{\lg{\vec\xt}}),
\end{equation}
where missing parentheses are to be added on
the right: $\xx \op \yy \op \zz$ stands for $\xx \op
(\yy \op \zz)$.
\end{nota}

\begin{lemm} \label{L:Seq}
Assume that $(\SS, \op)$ is an LD-system. Then
$(\Seq\SS, \ops, \concat)$ is an ALD-system.
\end{lemm}

\begin{proof}
The only point that is not absolutely obvious is that $(\LD)$
holds. Now, for all $\xs$, $\xt$, $\xu$ in~$\Seq\SS$, the
$k$th entry in~$\vec\xs \ops (\vec\xt \ops \vec\xu)$
is $\xs_1 \op \dd \op \xs_p \op \xt_1 \op \dd \op
\xt_q \op \xu_k$, while that of
$(\vec\xs \ops \vec\xt) \ops (\vec\xs \ops \vec\xu)$ is
$$(\xs_1 \op \dd \op \xs_p \op \xt_1) \op \dd \op
(\xs_1 \op \dd \op \xs_p \op \xt_q) \op \xs_1 \op \dd
\op \dd \op \xu_k.$$
Repeated applications of the LD law show that the
expressions are equal.
\end{proof}

We can now introduce the two mappings that give rise
to ALD-invariants.

\begin{defi}
For each term~$\tt$ in~$\TTOo$, we define a
term~$\II\tt$ in~$\TO1$ and a finite sequence of
terms~$\JJ\tt$ in~$\Seq\TTo$ using the inductive
clauses
\begin{gather}
(\II\tt, \JJ\tt) = \begin{cases}
(\xx, \tt)&\mbox{if $\tt$ is a variable,}\\
(\II\ttii, \JJ\tti \ops \JJ\ttii)&\mbox{for $\tt = \tti \op \ttii$,}\\
(\II\tti \OP \II\ttii, \JJ\tti \concat \JJ\ttii)
&\mbox{for $\tt = \tti \OP \ttii$,}
\end{cases}
\end{gather}
\end{defi}

For instance, for $\tt = \xx_1 \op ((\xx_2 \op \xx_3)
\OP \xx_4)$, the reader can check the values
$\II\tt =  \xx \OP \xx$, $\JJ\tt = (\xx_1 \op (\xx_2 \op
\xx_3), \xx_1 \op \xx_4)$.

\begin{lemm} \label{L:Invariant}
Assume that $\tt, \tt'$ are ALD-equivalent terms
in~$\TTOo$. Then we have
\begin{equation}
\II\tt = \II{\tt'}
\mbox{\quad and \quad}
\JJ\tt \eLDs \JJ{\tt'},
\end{equation}
the latter meaning that the sequences~$\JJ\ts$ and
$\JJ\tt$ have equal lengths and pairwise
$\LD$-equivalent entries.
\end{lemm}

\begin{proof}
As ALD-equivalence is the congruence on~$\TTOo$
generated by the pairs of terms occurring in the
laws~$(\LD)$, $(\ALD_1)$, and~$(\ALD_2)$, it is sufficient
to check that the relations $\II\tt = \II{\tt'}$ and $\JJ\tt
\eLD \JJ{\tt'}$ are congruences on~$\TTOo$, and that they
include all instances of~$(\LD)$, $(\ALD_1)$, and~$(\ALD_2)$.

The fact that~$\II{\tti \op \ttii}$ and~$\II{\tti \OP \ttii}$ are
defined from~$\II\tti$ and~$\II\ttii$ makes it clear that $\II\tt
= \II{\tt'}$ is a congruence, \ie, that it is compatible
with~$\op$ and~$\OP$.  The same argument works for
$\JJ\tt \eLD \JJ{\tt'}$, as the relation~$\eLD$
on~$\Seq\TTo$ is itself a congruence.

Let $(\tt, \tt')$ be an instance of~$(\LD)$, \ie, assume that $\tt$
and¨$\tt'$ are of the form
$\tt =  \tti \op (\ttii \op \ttiii)$ and $\tt' =  (\tti \op
\ttii) \op (\tti \op \ttiii)$. The
definitions yields
\begin{gather*}
\II\tt = \II\ttiii 
= \II{\tt'},\\
\JJ\tt = \JJ\tti \ops (\JJ\ttii \ops
\JJ\ttiii), \quad 
\JJ{\tt'} = 
(\JJ\tti \ops \JJ\ttii) \ops (\JJ\tti \ops
\JJ\ttiii),
\end{gather*}
and the latter are $\eLD$-equivalent by
Lemma~\ref{L:Seq}.
Similarly, for $(\tt, \tt')$ an instance of~$(\ALD_1)$,
\ie, for $\tt = \tti \op (\ttii \op \ttiii)$ and $\tt' =
(\tti \OP \ttii) \op \ttiii$, we have
$$\II\tt = \II\ttiii = \II{\tt'}
\mbox{\quad and \quad}
\JJ\tt = \JJ\tti \ops \JJ\ttii \ops \JJ\ttiii = \JJ{\tt'}.$$
Finally, for $(\tt, \tt')$ an instance of~$(\ALD_2)$, \ie,
for $\tt = \tti \op (\ttii \OP \ttiii)$ and $\tt' = (\tti \op
\ttii) \OP (\tti \op \ttiii)$, we find
$$\II\tt = \II\ttii \OP \II\ttiii = \II{\tt'}
\mbox{\quad and \quad}
\JJ\tt = (\JJ\tti \ops \JJ\ttii) \concat (\JJ\tti \ops
\JJ\ttiii) =
\JJ{\tt'},$$ which completes the proof.
\end{proof}

\begin{rema}
The result that $I$ is an ALD-invariant can also be
deduced from applying the construction of Example¨\ref{X:ALD}
to the free algebra $(\TO1, \op)$---as well as the result that $J$
mod.\,$(\LD)$ is an ALD-invariant follows from the construction
of Lemma¨\ref{L:Seq}.
\end{rema}

\subsection{Special terms}

We shall now see that, for each term~$\tt$ in~$\TTOo$,
the pair $(\II\tt, \JJ\tt)$ determines the ALD-class
of~$\tt$.

\begin{defi}
For $\tv$ is a term of size~$p$ in~$\TO1$, and
$\vtt$ is a length~$p$ sequence of terms
in~$\TTo$, we denote by $\tv[\vtt]$ the term obtained
from~$\tv$ by substituting~$\tt_1\ddd\tt_p$ to the variables
of~$\tt$ enumerated from left to right. A term is called {\it
special} if it is of the form¨$\tv[\vtt]$ with¨$\tv, \vtt$ as above.
\end{defi}

Saying that a term~$\tt$ is special means that, in the tree
associated with~$\tt$, no~$\OP$ symbol lies below an $\op$
symbol (according to the convention that the root lies on the
top).  The following result shows non only that every term
in~$\TTOo$ is ALD-equivalent to a special term, but also that the
pair $(\II\tt, \JJ\tt)$ determines the ALD-class of~$\tt$.

\begin{lemm} \label{L:Special}
For every term~$\tt$ in~$\TTOo$ we have
\begin{equation} \label{E:Special}
\tt \eALD \II\tt[\JJ\tt].
\end{equation}
\end{lemm}

\begin{proof}
If $\tt$ is a variable, \eqref{E:Special} is an equality. For an
induction, it is sufficient to show that the following
relations hold for all terms¨$\tu, \tv$ in¨$\TO1$ and
all sequences¨$\vts, \vtt$ in¨$\Seq\TTo$
\begin{gather}
\label{E:Special1}
\tu[\vts] \op \tv[\vtt] \eALD \tv[\vts \ops \vtt],\\
\label{E:Special2}
\tu[\vts] \OP \tv[\vtt] = (\tu \OP \tv)[\vts \concat \vtt].
\end{gather}

We establish \eqref{E:Special1} using  induction on the sum of the
sizes, say~$p$ and~$q$,  of~$\tu$ and~$\tv$, which also are the
lengths of¨$\vts$ and¨$\vtt$, respectively. We recall that
missing parentheses are to be added on the right, \ie, $\xx \op
\yy \op \zz$ stands for $(\xx \op \yy) \op (\xx \op \zz)$.

For $p = q = 1$, the terms~$\tu$ and~$\tv$ are variables, so we
have $\tu[\vts] = s_1$ and $\tv[\vtt] = t_1$, and
\eqref{E:Special1} reduces to the equality $s_1 \op t_1 =
\tv[s_1 \op t_1]$. Assume now $p +
q > 2$. Then we have $p \ge 2$ or $q \ge 2$. Assume
first $q \ge 2$. Write $\tv = \tv_1 \OP \tv_2$, and let
$r$ be the size of~$\tv_1$. We find
\begin{align*}
\tu[\vts] \op \tv[\vtt]
&= \tu[\vts]
\op (\tv_1[t_1\ddd t_r] \OP
\tv_2[t_{r+1}\ddd t_q]) 
&\mbox{(by definition)}\\
&\eALD
(\tu[\vts]
\op \tv_1[t_1\ddd t_{r}])  \OP (\tu[\vts]
\op \tv_2[t_{r+1}\ddd t_q])) 
&\mbox{$(\ALD_2)$}\\
&\eALD
\tv_1[\vts \ops t_1\ddd \vts \ops t_r] 
\OP \tv_2[\vts \ops 
t_{r+1}\ddd \vts \ops t_q] 
&\mbox{(by ind. hyp.)}\\
&= (\tv_1 \OP \tv_2))[\vts \ops \vtt]
&\mbox{(by definition).}
\end{align*}
Assume now $p \ge 2$. Writing similarly $\tu = \tu_1 \OP
\tu_2$, and  letting $r$ be now the size of~$\tu_1$, we find
\begin{align*}
\tu[\vts] \op \tv[\vtt]
&= (\tu_1[s_1\ddd s_r] \OP \tu_2[s_{r+1}\ddd s_p])
\op \tv[\vtt]
&\mbox{(by definition)}\\
&\eALD
\tu_1[s_1\ddd s_r] \op (\tu_2[s_{r+1}\ddd s_p] \op \tv[\vtt])
&\mbox{$(\ALD_1)$}\\
&\eALD
\tu_1[s_1\ddd s_r] \op
\tv[s_{r+1} \op\dd \op s_p\ops \vtt]
&\mbox{(by ind. hyp.)}\\
&\eALD
\tv[s_1 \op \dd \op s_p \op s_{r+1} \op\dd \op s_p\ops \vtt]
= \tv[\vts \ops \vtt]
&\mbox{(by ind. hyp.).}
\end{align*}
As for¨\eqref{E:Special2}, it follows from the definition directly. 
\end{proof}

\subsection{The word problem of~$\ALD$}

It is now easy to solve the word problem for¨$\ALD$.

\begin{prop} \label{P:Word}
The word problem of¨$\ALD$ is decidable: if $\tt, \tt'$ are
terms in¨$\TTOo$, then $\tt \eALD \tt'$ holds if and only if the
terms¨$\II\tt$ and¨$\II{\tt'}$ are equal, and the
sequences¨$\JJ\tt$ and¨$\JJ{\tt'}$ have the same length and
consist of pairwise LD-equivalent terms of¨$\TTo$.
\end{prop}

\begin{proof}
The condition is necessary by Lemma¨\ref{L:Invariant}. It is
sufficient by Lemma¨\ref{L:Special}. Indeed, if $\vts$, $\vtt$ are
length¨$p$ sequences of pairwise LD-equivalent terms in¨$\TTo$
and if $\tv$ is any size¨$p$ term in¨$\TO1$, the
terms¨$\tv[\vts]$ and¨$\tv[\vtt]$ are ALD-equivalent. So, if
$\tt, \tt'$ are terms in¨$\TTOo$ satisfying $\II\tt = \II{\tt'}$ and
$\JJ\tt \eLD \JJ{\tt'}$, we obtain $$\tt \eALD \II\tt[\JJ\tt] \eALD
\II{\tt'}[\JJ{\tt'}] \eALD \tt',$$ hence $\tt \eALD \tt'$.
As the relation¨$\eLD$ is known to be decidable \cite{Dgd}, so
is¨$\eALD$.
\end{proof}

As for the complexity of the previous solution, the known
upper bounds for the word problem of¨$(\LD)$ are a single
exponential in the case of terms with one variable, and a double
exponential in the general case. As the size of the
sequence¨$\JJ\tt$ may be exponential in the length of¨$\tt$
since each application of¨$(\ALD_2)$ may double the
length, the solution described in Proposition¨\ref{P:Word} has a
(certainly not optimal) upper bound which is doubly
exponential in the case of one variable, and triply exponential in
the general case---the results of Section¨\ref{S:Braids} below will
give a better, simply exponential algorithm in the case of one
variable.

\section{Parenthesized braids} \label{S:Braids}

The group of parenthesized braids~$\Bb$ was
introduced in~\cite{Bri1, Bri2, Dhb}---in a
different framework---and further investigated
in~\cite{Dhe}. It is shown in the latter paper that
$\Bb$ can be equipped with two binary operations
that make it an ALD-system. The aim of this section is to study
this specific ALD-system, and in particular to show that it
contains many copies of the free ALD-system on one generator.

\subsection{The group~$\Bb$}

The simplest way to introduce¨$\Bb$ is to start from a
presentation:

\begin{defi}
We denote by¨$\Bb$ the group generated by two infinite
sequences¨$\ss1, \ss2, \dd$, $\aa1, \aa2, \dd$ subject to the
relations
\begin{equation} \label{E:Present}
\begin{cases}
\quad \ss j \ss i = \ss i \ss j  , \quad \aa j \ss i = \ss i \aa j 
\qquad
&\mbox{for $j \ge i + 2$},\\
\quad \aa j \ss i = \ss{i+1} \aa j , \quad 
\aa j \aa j = \aa{i+1} \aa j \quad
&\mbox{for $j \le  i - 1$},\\
\quad \ss j \ss i \ss j = \ss i \ss j \ss i, \quad
\ss i \ss j \aa i = \aa j \ss i, \quad
\ss j \ss i \aa j = \aa i \ss i, \quad
&\mbox{for $j = i+1$}.
\end{cases}
\end{equation}
\end{defi}

\begin{figure} [htb]
\begin{picture}(83,22)(0,0)
\put(10,10){\includegraphics{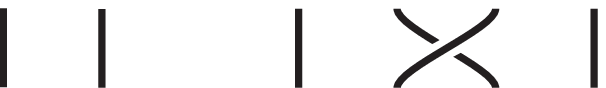}}
\put(-4,13){$\ss i\mapsto$}
\put(9.5,19){$1$}
\put(19.5,19){$2$}
\put(50,19){$i$}
\put(58,19){$i{+}1$}
\put(27,13){\dots}
\put(77,13){\dots}
\put(10,0){\includegraphics{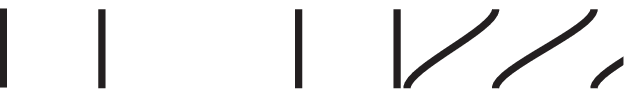}}
\put(-4,3){$\aa i\mapsto$}
\put(27,3){\dots}
\put(77,3){\dots}
\end{picture}
\caption{\smaller Diagram representation of the generators
of¨$\Bb$: there are infinitely strands numbered by positive
integer coefficients polynomials in an infinitely small
variable¨$\epsilon$; the effect of¨$\ss i$ is to let all strands
with index $i + 1 + o(1)$ cross over all strands with
index¨$i + o(1)$; the effect of¨$\aa i$ is to shrink
all strands of the form $i + o(1)$ by a factor¨$\epsilon$ and to
left translate all strands with index¨$\ge i+1$ so as to avoid
gaps.}
\label{F:Diag}
\end{figure}

It is shown in¨\cite{Bri2} that $\Bb$ is actually generated
by¨$\ss1, \ss2, \aa1$, and¨$\aa2$, and that it admits a
finite---but much less readable---presentation with respect to
those generators.  It is shown in¨\cite{Dhe} that the elements
of¨$\Bb$ admit a natural geometric interpretation in terms of
parenthesized braid diagrams, which are similar to ordinary braid
diagrams---{\it cf.}¨for instance \cite{Bir, Dgd, PrS}---but with
non-uniform distances between the strands. As we shall use this
interpretation here---nor do wo either use the interpretation in
terms of isotopy classes of homeomorphisms of a sphere with
a Cantor set of punctures---we shall not go into details here and
just refer to Figure¨\ref{F:Diag} for a rough intuition.

\begin{defi}
We denote by¨$\sh$ the endomorphism of the group¨$\Bb$
that maps $\ss i$ to¨$\ss{i+1}$ and $\aa i$ to¨$\aa{i+1}$ for
each¨$i$.
\end{defi}

It is shown in¨\cite{Dhe} that $\sh$ is injective---but not
surjective: neither $\ss1$ nor¨$\aa1$ belong to¨$\Im\sh$.

\begin{figure} [htb]
\begin{picture}(108,39)(0,0)
\put(0,0){\includegraphics{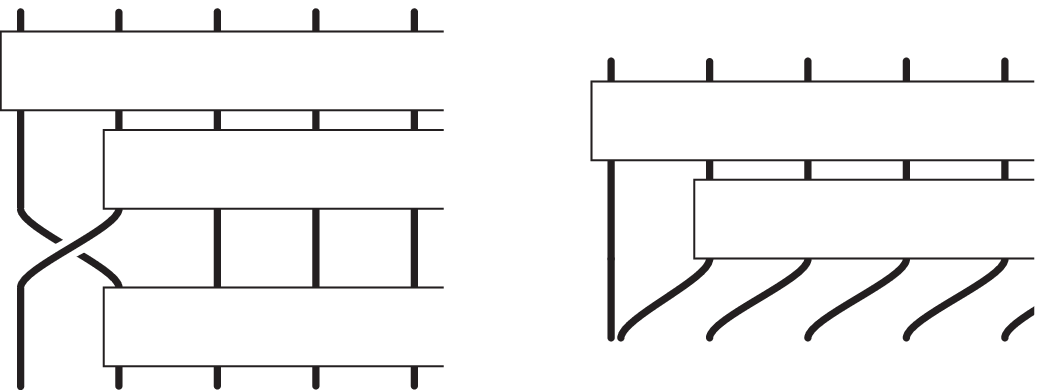}}
\put(21,32){$\bx$}
\put(26,22){$\by$}
\put(25,5){$\bx^{-1}$}
\put(81,26.5){$\bx$}
\put(86,17){$\by$}
\put(44,6){$\dots$}
\put(44,23){$\dots$}
\put(44,33){$\dots$}
\put(104,17){$\dots$}
\put(104,27){$\dots$}
\end{picture}
\caption{\smaller Diagram representation of the ALD
operations on¨$\Bb$: the diagram of $\bx \op \by$ (left) and
$\bx \OP \by$ (right) from those of¨$\bx$ and¨$\by$}
\label{F:Oper}
\end{figure}

\begin{prop} [\cite{Dhe}] \label{P:Braids}
(Figure¨\ref{F:Oper}) Let $\op$, $\OP$ be the binary
operations on¨$\Bb$ defined by
\begin{equation} \label{E:Braids}
\bx \op \by := \bx \oop \sh\by \oop \ss1 \oop \sh\bx\inv,
\qquad
\bx \OP \by := \bx \oop \sh\by \oop \aa1.
\end{equation}
Then $(\Bb, \op, \OP)$ is an ALD-system.
\end{prop}

\subsection{A freeness criterion}

Our aim is to show that the ALD-system $(\Bb, \op, \OP)$
includes copies of the free ALD-system of rank¨1. To prove the
result, we need a criterion for recognizing such free ALD-systems.

Assume that $(\SS, \op, \OP)$ is a double binary system
generated by a single element¨$\gg$. Then, there exists a
surjective homomorphism¨$\pi$ of ¨$\TOo1$ onto¨$\SS$ that
maps¨$\xx$ to¨$\gg$: by definition, the value¨$\pi(\tt)$ is the
{\it evaluation} of¨$\tt$ at¨$\gg$, and it will be denoted
by¨$\tt(\gg)$---exactly as the evaluation of a polynomial¨$P$
at¨$\gg$ would be denoted by¨$P(\gg)$. Then, saying that
$(\SS, \op, \OP)$ is an ALD-system means that $\tt \eALD \tt'$
implies $\tt(\gg) = \tt'(\gg)$, and saying that $(\SS, \op, \OP)$
is a free ALD-system based on¨$\{\gg\}$ means that $\tt \eALD
\tt'$ is equivalent to $\tt(\gg) = \tt'(\gg)$. In other words, in
roder to prove that some ALD-system¨$\SS$ generated by an
element¨$\gg$ is free, the point is to prove that $\tt(\gg)  \not=
\tt'(\gg)$ holds for all pairs of terms¨$(\tt, \tt')$ satisfying $\tt
\not\eALD \tt'$. The criterion we shall establish new allows one to
restrict to pairs of terms¨$(\tt, \tt')$ of a restricted type.

\begin{defi}
For $\tu, \tv$ in~$\TO1$, we say that $\tu \sm
\tv$ holds if we have either $(i)$ $\tu = \xx$ and
$\tv \not= \xx$,  or $(ii)$ $\tu = \tu_1 \OP \tu_2$
and $\tv = \tv_1 \OP \tv_2$ with $\tu_1 \sm
\tv_1$, or $(iii)$ $\tu = \tu_1 \OP \tu_2$
and $\tv = \tv_1 \OP \tv_2$ with $\tu_1 =
\tv_1$ and $\tu_2 \sm \tv_2$.
\end{defi}

Clearly, the relation~$\sm$ is a strict linear order
on~$\TO1$.

\begin{defi}
$(i)$ For $\ts, \tt$ in~$\TTo$, we say that $\ts
\pref \tt$ holds if there exist $p \ge 1$ and terms
$\tt_1 \ddd \tt_p$ in~$\TTo$ satisfying 
$$\tt = (\dd((\ts \op \tt_1) \op \tt_2) \dd ) \op
\tt_p.$$

$(ii)$ For $\vts, \vtt$ in~$\Seq\TTo$,  we say that
$\vts \prefs \vtt$ holds if the lengths of¨$\vts$ and¨$\vtt$
are equal and there exists¨$k \le \lg\vts$ satisfying $\ts_i =
\tt_i$ for~$i < k$ and $\ts_k \pref \tt_k$.
\end{defi}

\begin{prop}\label{P:Criterion}
Assume that $\SS$ is an ALD-system generated by
an element~$\gg$. Then a necessary and sufficient condition for
$\SS$ to be free based on~$\{\gg\}$ is that $\SS$ satisfies no
equality of the form
\begin{equation} \label{E:Criterion}
\tu[\vts](\gg) = \tv[\vtt] (\gg)
\end{equation}
with $\tu, \tv$ in¨$\TO1$ and $\vts, \vtt$ in¨$\Seq{\To1}$
satisfying either $\tu \sm \tv$, or $\tu = \tv$ and $\vec\ts
\prefs \vtt$.
\end{prop}

\begin{proof}
Assume that $\ts, \tt$ are $\ALD$-inequivalent terms
in~$\TOo1$. As was said above, the problem is to show that the
evaluations $\ts(\gg)$ and $\tt(\gg)$ of~$\ts$ and~$\tt$
in~$\SS$ cannot be equal. By Lemma~\ref{L:Special}, there exist
$\tu, \tv$ in~$\TO1$ and $\vec\ts$, $\vtt$ in~$\Seq\TTo$
satisfying $\ts \eALD \tu[\vts]$ and $\tt \eALD
\tv[\vtt]$. As $\SS$ is an ALD-system, we have $\ts(\gg) =
\tu[\vts](\gg)$ and $\tt(\gg) = \tv[\vtt](\gg)$,
so it is sufficient to prove $\tu[\vts](\gg) \not=
\tv[\vtt] (\gg)$. Now, by Lemma~\ref{L:Free}, the
hypothesis $\tu[\vts] \not\eALD \tv[\vtt]$
implies $\tu \not= \tv$, or $\tu = \tv$ and $\vec\ts
\not\eLDs \vtt$. In the first case, we must have either $\tu \sm
\tv$ or $\tv \sm \tu$ as $\sm$ is a linear ordering, hence, if
no equality¨\eqref{E:Criterion} holds, we deduce $\ts(\gg) \not=
\tt(\gg)$. In the second case, as the sequences
$\vec\ts$ and $\vtt$ have the same length, there exists an
index~$k \le p$ such that we have $\ts_i \eLD \tt_i$ for $i < k$
and $\ts_k \not\eLD \tt_k$. By the results of~\cite{Dgd}, the
latter relation implies the existence of terms~$\ts'_k, \tt'_k$
satisfying $\ts'_k \eLD \ts_k$, $\tt'_k \eLD \tt_k$ and either
$\ts'_k \pref \tt_k'$ or $\tt'_k \pref \ts'_k$. Let~$\vec\ts'$
denote the sequence obtained from~$\vec\ts$ by
replacing~$\ts_k$ by~$\ts'_k$, and let~$\vtt'$ denote the
sequence obtained from~$\vtt$ by replacing~$\tt_i$
with~$\ts_i$ for~$i < k$, and by replacing~$\tt_k$
with~$\tt'_k$. Then, as $\SS$ is an ALD-system, we have 
$\tu[\vts](\gg) = \tu[\vts'](\gg)$ and
$\tu[\vtt](\gg) = \tu[\vtt'](\gg)$, and, by
construction, we have $\vec\ts' \prefs \vtt'$ or $\vtt'
\prefs \vec\ts'$. If no equality¨\eqref{E:Criterion} holds, we
deduce $\tu[\vts'](\gg) \not= \tu[\vtt'](\gg)$, hence $\ts(\gg)
\not= \tt(\gg)$.
\end{proof}

\subsection{Term evaluation}

In order to apply the criterion of Proposition¨\ref{P:Criterion} in
the ALD-system $(\Bb, \op, \OP)$, we need to be able to evaluate
in¨$\Bb$ expressions of the form¨$\tv[\vtt](\gg)$ with¨$\tv$
a term in¨$\TO1$ and $\vtt$ a sequence of terms in¨$\To1$.
To this end, we shall use the following explicit formulas.

\begin{lemm} \label{L:Shift}
Assume that $\tv$ is a term of size~$p$ in~$\TO1$. Then, for
each~$\bx$ in~$\Bb$, we have
\begin{equation} \label{E:Shift}
\tv(1) \oop \sh\bx = \sh^p\bx \oop \tv(1).
\end{equation}
\end{lemm}

\begin{proof}
We use induction on~$\tv$. For $\tv = \xx$, we have $p = 1$ and
$\tv(1) = 1$, so \eqref{E:Shift} is true. Otherwise, assume
$\tv = \tv_1 \OP \tv_2$. By definition, 
$\tv(1)$ is $\tv_1(1) \oop \sh{\tv_2(1)}
\oop \aa1$. Let
$p_i$ be the size of~$\tv_i$. Using the induction
hypothesis, we find
\begin{align*}
\tv(1) \oop \sh{\bx} 
&=  \tv_1(1) \oop \sh{\tv_2(1)} \oop \aa1
\oop \sh{\bx}
&\mbox{(by definition)} \\
&=  \tv_1(1) \oop \sh{\tv_2(1)}  \oop
\sh^2{\bx} \oop \aa1
&\mbox{(by the relations of¨$\Bb$)}\\
&=  \tv_1(1) \oop \sh({\tv_2(1)}  \oop
\sh{\bx}) \oop \aa1\\  
&=  \tv_1(1)   \oop \sh(\sh^{p_2}{\bx} \oop
\tv_2(1)) \oop \aa1
&\mbox{(by induction hypothesis)}\\ 
&= \tv_1(1)  \oop \sh(\sh^{p_2}{\bx}) \oop
\sh{\tv_2(1)} \oop \aa1 \\
&= \sh^{p_1}{(\sh^{p_2}{\bx})} \oop
\tv_1(1)   \oop
 \sh{\tv_2(1)} \oop \aa1
&\mbox{(by induction hypothesis)}
\end{align*}
and the latter is~$\sh^p{\bx} \oop \tv(1)$.
\end{proof}

\begin{lemm} \label{L:Eval}
Assume $\tt = \tv[\vtt]$, with $\tv$ a size¨$p$ term
in¨$\TO1$ and $\vtt$ a length¨$p$ sequence of terms
in¨$\To1$. Then, for each¨$\gg$ in¨$\Bb$, we have
\begin{equation} \label{E:Eval}
\tt(\gg) = \tt_1(\gg) \oop \sh\tt_2(\gg) \oop \dd 
\oop \sh^{p-1}\tt_p(\gg) \oop \tv(1).
\end{equation}
\end{lemm}

\begin{proof}
We use induction on~$\tv$. For $\tv = \xx$,
we have $p = 1$ and $\tt = \tt_1$, so the result is
clear. Otherwise, assume $\tv = \tv_1 \OP
\tv_2$. Let $q$ be the size of~$\tv_1$. Then 
we have 
$$\tt = \tv_1[\tt_1 \ddd \tt_q] \OP
\tv_2[\tt_{q+1} \ddd \tt_p],$$
and, using the induction hypothesis twice,
we deduce
\begin{align*}
\tt(\gg) 
&=  \tv_1[\tt_1\ddd \tt_q](\gg)
\oop \sh\tv_2[\tt_{q+1}\ddd \tt_p](\gg)
\oop \aa1 \\ 
&=\tt_1(\gg) \oop \dd 
\oop \sh^{q-1}\tt_q(\gg) \oop \tv_1(1) 
\oop \sh\tv_2[\tt_{q+1}\ddd \tt_p](\gg)
\oop \aa1 
&\mbox{(ind. hyp.)}\\
&=\tt_1(\gg) \oop \dd 
\oop \sh^{q-1}\tt_q(\gg)  
\oop \sh^q\tv_2[\tt_{q+1}\ddd \tt_p](\gg)
\oop \tv_1(1) \oop \aa1 
&\mbox{\eqref{E:Shift}}\\
&=
 \tt_1(\gg) \oop \dd 
\oop \sh^{q-1}\tt_q(\gg)  
\oop \sh^q\tt_{q+1}(\gg) \oop \dd
\oop \sh^{p-1}\tt_p(\gg)  \oop \sh^q\tv_2(1)
\oop \tv_1(1) \oop \aa1 
&\mbox{(ind. hyp.)}\\
&= \tt_1(\gg) \oop \dd 
\oop \sh^{q-1}\tt_q(\gg)  
\oop \sh^q\tt_{q+1}(\gg) \oop \dd
\oop \sh^{p-1}\tt_p(\gg)
\oop \tv_1(1)   \oop \sh\tv_2(1) \oop
\aa1,
&\mbox{\eqref{E:Shift}}
\end{align*}
which gives~\eqref{E:Eval} since we have $\tv(1) = \tv_1(1) \oop
\sh\tv_2(1) \oop \aa1$.
\end{proof}

\subsection{Monogenerated subsystems of~$\Bb$}

It is shown in~\cite{Dhe} that the evaluation mapping $\tv
\mapsto \tv(1)$ of~$\TO1$ into~$\Bb$ is injective. We shall need
the following strengthening of this result:

\begin{lemm} \label{L:Image}
If $\tu, \tv$ are distinct terms in~$\TO1$, then, in~$\Bb$, the
quotient $\tu(1)\inv \tv(1)$ does not belong to~$\Im\sh$.
\end{lemm}

\begin{proof}
Let $\xx^{[N]}$ denote the term
of~$\TO1$ inductively defined by $\xx^{[1]} = \xx$ and
$\xx^{[N]} = \xx \OP \xx^{[N-1]}$ for $N \ge 2$. The subgroup
of¨$\Bb$ generated bu the elements¨$\aa i$ is isomorphic to
Thompson's group¨$F$, and it gives rise to a partial action
on¨$\TO1$ corresponding to applying the associativity
law~\cite{Dhe}: the action of¨$\aa i$ on a term¨$\tv$ is defined
provided¨$\tv$ can be expressed as $\tv_1 \OP \dd \OP
\tv_{i+2}$, \ie, we have $\hR(\tv) \ge i+2$, and, in this case, one
defines $\tv \cdot \aa i = \tv_1 \OP \dd \OP \tv_{i-1} \OP (\tv_i
\OP \tv_{i+1}) \OP \tv_{i+2}$. Then, an easy induction shows
that, for each term~$\tv$ of size~$p$ in~$\TO1$, the
element¨$\tv(1)$ of¨$\Bb$ maps any sufficiently large
term¨$x^{[N]}$ to the term~$\tv \OP x^{[N-p]}$. Hence
$\tu(1)\inv \tv(1)$ maps $\tu \OP \xx^{[N -p]}$ to $\tv \OP
\xx^{[N-q]}$, where $p$ is the size of~$\tu$. Now any element
of~$\Im\sh$ maps a term of the form $\tu \OP \dd$ to
another term of the form~$\tu \OP \dd$, since only $\aa1$ may
change the left subterm of the initial term. Hence $\tu(1)\inv
\tv(1) \in \Im\sh$ is impossible for $\tv \not=
\tu$.
\end{proof}

\begin{prop} \label{P:Main}
For any~$\gg$ in~$\Bb$, the
closure of~$\{\gg\}$ under~$\op$ and~$\OP$ is free
ALD-system.
\end{prop}

\begin{proof}
We apply the criterion of Proposition~\ref{P:Criterion}. Assume
that $\tu, \tv$ are terms in~$\TO1$ and $\vts, \vtt$
are sequences of terms in~$\To1$. Let $\ts = \tu[\vts]$
and $\tt = \tv[\vtt]$.  Our aim is to prove $\ts(\gg)\inv
 \tt(\gg) \not= 1$ both for $\tu \sm \tv$, and for $\tu =
\tv$ with $\vec\ts \prefs \vtt$. Applying
Lemma~\ref{L:Eval}, we find
\begin{equation} \label{E:Comp}
\ts(\gg)\inv \tt(\gg) = 
\tu(1)\inv \oop \sh^{p-1}\ts_p(\gg)\inv \oop\dd\oop 
\ts_1(\gg)\inv \oop 
\tt_1(\gg) \oop\dd\oop \sh^{q-1}\tt_q(\gg) \oop \tv(1).
\end{equation}
We shall consider three cases, which cover the cases
$\tu \sm \tv$, and $\tu = \tv$ with $\vec\ts \prefs \vtt$,
and prove in each of them that the right hand side
of¨\eqref{E:Comp} is not¨$1$.

Assume first that there exists $k \le \inf(p, q)$ such that $\ts_i
\eLD \tt_i$ holds for $i < k$, and $\ts_k \not\eLD \tt_k$ holds.
Then we have $\ts_i(\gg) = \tt_i(\gg)$ for $i < k$, and
\eqref{E:Comp} becomes
$$\ts(\gg)\inv \tt(\gg) = 
\tu(1)\inv \oop \sh^{p-1}\ts_p(\gg)\inv \oop\dd\oop 
\sh^{k-1}(\ts_k(\gg)\inv \tt_k(\gg)) \oop\dd\oop
\sh^{q-1}\tt_q(\gg) \oop
\tv(1).$$
By the results of¨\cite{Dgd}, the hypothesis $\ts_k \not\eLD
\tt_k$ implies either $\ts_k \pLD \tt_k$ or $\tt_k \pLD \ts_k$,
and the explicit definition of operation¨$\op$ on¨$\Bb$ then
implies that the braid $\ts_k(\gg)\inv \tt_k(\gg)$ admits an
expression where the generator¨$\ss1$ appears but¨$\ss1\inv$
does not, or $\ss1\inv$ appears but
$\ss1$ does not. It follows that $\ts(\gg)\inv \tt(\gg)$ admits an
expression in which $\ss k$ appears but neither¨$\ss k\inv$ nor
any¨$\ss i^{\pm1}$ with $i < k$ does, or vice versa exchanging
$\ss k$ and¨$\ss k\inv$. By¨\cite{Dhe}, Proposition¨4.6, this
guarantees $\ts(\gg) < \tt(\gg)$ in the canonical ordering
of¨$\Bb$, hence $\ts(\gg) \not= \tt(\gg)$.

Assume now $p < q$ with $\ts_i \eLD \tt_i$ for $i \le p$. In this
case, \eqref{E:Comp} reduces to
\begin{equation*} 
\ts(\gg)\inv  \tt(\gg) = 
\tu(1)\inv \oop \sh^p(\tt_{p+1}(\gg) \oop\dd\oop
\sh^{q-p-1}\tt_q(\gg)) \oop \tv(1).
\end{equation*}
By Lemma~\ref{L:Shift}, we have $\tu(1)\inv \oop \sh^p\zz =
\sh\zz \tu(1)\inv$ for each~$\zz$ in~$\Bb$, so we get
$$\ts(\gg)\inv \oop \tt(\gg) =  
\sh(\tt_{p+1}(\gg) \oop\dd\oop
\sh^{q-p-1}\tt_q(\gg)) \oop \tu(1)\inv\tv(1).$$
This cannot be~$1$, as the first factor belongs to~$\Im\sh$,
while, according to Lemma~\ref{L:Image}, $\tu(1)\inv\tv(1)$
does not unless $\tu = \tv$ holds.

Assume finally $p = q$ with $\ts_i \eLD \tt_i$ for $i \le p$, and
$\tu \sm \tv$. Then \eqref{E:Comp} reduces to
\begin{equation*} 
\ts(\gg)\inv \tt(\gg) = 
\tu(1)\inv \oop \tv(1),
\end{equation*}
and, by Lemma~\ref{L:Shift}, the above expression cannot
be¨$1$.
\end{proof}

\begin{rema}
It is shown in¨\cite{Dhe} that the parenthesized braid
group~$\Bb$ comes can be equipped with a distinguished linear
ordering that extends both the linear ordering of braids
and the natural ordering on Thompson's group induced by the
lexicographical ordering of finite trees. Let us define a
relation¨$\smALD$ on special terms in¨$\TOo1$ as follows:
first say that $\tu[\vts] \sm \tv[\vtt]$ holds if we have either
$\vts \prefs \vtt$, or $\vts$ is a proper prefix of¨$\vtt$, or
we have $\vts = \vtt$ and $\tu \sm \tv$ holds; then say that
$\ts \smALD \tt$ holds if there exist special terms
$\tu[\vts]$, $\tv[\vtt]$ satisfying $\tu[\vts] \sm \tv[\vtt]$,
$\ts \eALD \tu[\vts]$, and $\tt \eALD \tv[\vtt]$. Then
the relation¨$\smALD$ induces a linear ordering on the free
ALD-system $\TOo1/\!\!\eALD$, and what actually shows the
proof of Proposition¨\ref{P:Main} is that, for each parenthesized
braid¨$\gg$ in¨$\Bb$, the evaluation mapping
$\tt \mapsto \tt(\gg)$ is increasing.
\end{rema}

\subsection{The converse direction}

According to Proposition¨\ref{P:Braids}, the operations
of¨\eqref{E:Braids} define operations on¨$\Bb$ that make it an
ALD-system. We conclude with the easy observation that,
conversely, the operations defined on a group¨$G$ by formulas
of the type¨\eqref{E:Braids} give rise to an
ALD-system only if $G$ is closely connected to¨$\Bb$:

\begin{prop}
Assume that $G$ is a group, $\ssh$ is an endomorphism of¨$G$,
and $a, \sigma$ are fixed elements of¨$G$. Write $\ss i$ for
$\ssh^{i-1}(\s)$ and $\aa i$ for $\ssh^{i-1}(a)$. Then defining
\begin{equation} \label{E:OperBis}
\ex \op \ey = \ex \oop \ssh\ey \oop \s \oop \ssh\ex\inv,
\qquad 
\ex \OP \ey = \ex \oop \ssh\ey \oop a
\end{equation}
yields an ALD-system on the subgroup¨$H$ generated
by the elements¨$\ss i$'s and the¨$\aa i$'s---\ie, on the smallest
subgroup of¨$G$ containing¨$\s$ and¨$a$ and closed
under¨$\ssh$---if and only if the elements¨$\ss i$ and¨$\aa i$
obey the relations¨\eqref{E:Present}, \ie, if and only if $H$ is a
homomorphic image of¨$\Bb$.
\end{prop}

\begin{proof}
Assume that $(G, \op, \OP)$ is an ALD-system. The instance $1
\op (1 \op \ez) = (1 \op 1) \op (1 \op \ez)$ of¨$(\LD)$ expands
into
\begin{equation} \label{E:Commut}
\ssh^2\ez \oop \ss2\ss1 = \ss1 \oop \ssh^2\ez \oop
\ss2\ss1\ss2\inv.
\end{equation}
For $\ez = 1$, we obtain the braid relation
\begin{equation} \label{E:Rel1}
\ss1 \ss2 \ss1 = \ss2 \ss1 \ss2,
\end{equation}
and, then, \eqref{E:Commut} gives
\begin{equation}
\ssh^2\ez \oop \ss1 = \ss1 \oop \ssh^2\ez
\end{equation}
for each¨$\ez$. Similarly, the instance $1 \op (1 \op \ez) = (1
\OP 1) \op \ez$ of¨$(\ALD_1)$ expands into 
\begin{equation} \label{E:Commut2}
\ssh^2\ez \oop \ss2\ss1 = \aa1 \oop \ssh\ez \oop \ss1
\aa2\inv.
\end{equation}
For $\ez = 1$, we deduce 
\begin{equation} \label{E:Rel2}
\aa1\ss1 = \ss2\ss1\aa2,
\end{equation}
and, then, \eqref{E:Commut2} gives
\begin{equation}
\ssh^2\ez \oop \aa1 = \aa1 \oop \ssh\ez
\end{equation}
for each¨$\ez$. Finally, the instance $1 \op (1 \OP 1) = (1 \op
1) \OP (1 \op 1)$ of¨$(\ALD_2)$ expands into
\begin{equation} \label{E:Rel3}
\aa2\ss1 = \ss1\ss2\aa1.
\end{equation}
Conversely, it is easy to verify that the conjunction
of¨\eqref{E:Commut}, \eqref{E:Commut2} (for each¨$\ez$), 
and \eqref{E:Rel1}, \eqref{E:Rel2}, and \eqref{E:Rel3} guarantees
that $(G, \op, \OP)$ be an ALD-system. When we restrict to the
subgroup¨$H$, this amounts to saying that the elements¨$\ss i$
and¨$\aa i$ satisfy the defining relations¨\eqref{E:Present}
of¨$\Bb$.
\end{proof}

The previous result shows that there is no flexibility or
randomness in the construction of an ALD-system using the
formulas of¨\eqref{E:OperBis}. However, what was not
explained here---nor was it in¨\cite{Dhe} either---is where do
these formulas come from. Actually, the group¨$\Bb$ and the
formulas¨\eqref{E:OperBis} arise naturally when investigating
the so-called geometry monoid of the ALD laws. This will be
explained in a forthcoming paper.

\end{document}